\documentclass[11pt]{article}
\usepackage{amsfonts}
\usepackage{amsmath}
\usepackage{amssymb}
\usepackage{graphicx}

\newtheorem{theorem}{Theorem}[section]
\newtheorem{remark}{Remark}[section]
\newtheorem{definition}{Definition}[section]

\input{tcilatex}

\begin{document}

\author{Mohammed AL Horani\thanks{
Department of Mathematics, The University of Jordan, Amman, Jordan,
horani@ju.edu.jo. Current Address: Department of Mathematics, Faculty of
Science, University of Hail, Saudi Arabia, m.alhorani@uoh.edu.sa}, $\,$$\,$$\,$
Roshdi Khalil \thanks{
Department of Mathematics, The University of Jordan, Amman, Jordan,
roshdi@ju.edu.jo } and
Thabet Abdeljawad\thanks{
Department of Mathematics and Physical Sciences, Prince Sultan University, P.O.Box 66833, Riyadh 11586, Kingdom of Saudi Arabia,
tabdeljawad@psu.edu.sa}}
\title{Conformable Fractional Semigroups of Operators}
\date{}
\maketitle

\begin{abstract}
Let $X$ be a Banach space, and $T:[0,\infty )\rightarrow {\mathcal{L}}(X,X),$
the bounded linear operators on $X.$ A family $\{T(t)\}_{t\ge 0}\subseteq {%
\mathcal{L}}(X,X)$ is called a one-parameter semigroup if $T(s+t)=T(s)T(t),$
and $T(0)=I,$ the identity operator on $X.$ The infinitesimal generator of
the semigroup is the derivative of the semigroup at $t=0.$ The object of
this paper is to introduce a (conformable) fractional semigroup of operators whose
generator will be the fractional derivative of the semigroup at $t=0.$ The
basic properties of such semigroups will be studied.
\end{abstract}

\newline
\textbf{AMS classification numbers:} 26A33 \setcounter{MaxMatrixCols}{10}

\section{Introduction}

\label{intro} Let $X$ \ be a Banach space, and ${\mathcal{L}}(X,X)$ be the
space of all bounded linear operators on $X$. A family $\{T(t)\}_{t\ge
0}\subseteq {\mathcal{L}}(X,X)$ is called a one-parameter semigroup if:

\ \ \ \ \ \ \ \ \ \ \ \ \ \ \ \ \ \ \ \ \ \ \ \ \ \ \ \ $(i)\, T(0)=I,$ the
identity operator on $X.$

\ \ \ \ \ \ \ \ \ \ \ \ \ \ \ \ \ \ \ \ \ \ \ \ \ \ \ $\ (ii)\,
T(s+t)=T(s)T(t)\;\,$ for all $s,\; t\ge 0$. \newline
If, in addition, for each fixed $x\in X$, $T(t)x\longrightarrow x\;$ as $\;
t\longrightarrow 0^{+}$, then the semigroup is called c$_{0}-$semigroup or
strongly continuous semigroup. \newline

Semigroups of operators proved to be a very fruitful tool to solve
differential equations. One of the classical vector valued differential
equations is the so called the abstract Cauchy problem, precisely,
\begin{eqnarray*}
&&u(t)=Au(t),\quad t>0, \\
&&u(0)=x,
\end{eqnarray*}%
where $A:D(A)\subseteq X\longrightarrow X$ a linear operator of an
appropriate type, $x\in X$ is given and $\;u:[0,\infty )\longrightarrow X\;$
is the unknown function. We refer to \cite{EN} and \cite{Paz} for basic
theory of semigroups of operators and the abstract Cauchy problem. For the
inverse form of the abstract Cauchy problem, we refer to \cite{POV}, see
also \cite{H2010}.\newline

Fractional semigroups are related to the problem of fractional powers of
operators initiated first by Bochner, see \cite{Boch}. \ Balakrishnan, see
\cite{Bala}, studied the problem of fractional powers of closed operators
and the semigroups generated by them. \ The fractional Cauchy problem
associated with a Feller semigroup was studied by Popescu, see \cite{Pop}.%
\newline

In the literature, there are many definitions of fractional derivative. To
mention some:

\textbf{(i)} Riemann - Liouville Definition. For $\alpha \in \lbrack n-1,n),$
the $\alpha $ derivative of $f$ \ is
\begin{equation*}
{\large D}_{a}^{\alpha }{\large (f)(t)=}\frac{1}{{\huge \Gamma (n-\alpha )}}%
\frac{{\huge d}^{{\huge n}}}{{\huge dt}^{{\huge n}}}\underset{a}{\overset{t}{%
\int }}\frac{{\huge f(x)}}{{\huge (t-x)}^{{\huge \alpha -n+}1}}\,dx.
\end{equation*}

\textbf{(ii)} Caputo Definition. For $\alpha \in \lbrack n-1,n),$ the $%
\alpha $ derivative of $f$ \ is
\begin{equation*}
{\large D}_{a}^{\alpha }{\large (f)(t)=}\frac{1}{{\huge \Gamma (n-\alpha )}}%
\underset{a}{\overset{t}{\int }}\frac{{\huge f}^{{\large (n)}}{\huge (x)}}{%
{\huge (t-x)}^{{\huge \alpha -n+}1}}\,dx.
\end{equation*}

However, the following are the setbacks of one definition or the other:%
\newline
\newline
\textbf{(i)} The Riemann-Liouville derivative does not satisfy ${\large D}%
_{a}^{\alpha }{\large (1)=0}$ $({\large D}_{a}^{\alpha }{\large (1)=0}$ for
the Caputo derivative$)$, if $\alpha $ is not a natural number.\newline
\newline
\textbf{(ii)} All fractional derivatives do not satisfy the known formula of
the derivative of the product of two functions:
\begin{equation*}
{\large D}_{a}^{\alpha }{\large (fg)=fD}_{a}^{\alpha }{\large (g)+gD}%
_{a}^{\alpha }{\large (f).}
\end{equation*}%
\textbf{(iii)} All fractional derivatives do not satisfy the known formula
of the derivative of the quotient of two functions:
\begin{equation*}
\displaystyle{\large D}_{a}^{\alpha }{\large (f/g)=}\frac{{\large gD}%
_{a}^{\alpha }{\large (f)-fD}_{a}^{\alpha }{\large (g)}}{g^{2}}.
\end{equation*}%
\textbf{(iv)} All fractional derivatives do not satisfy the chain rule:
\begin{equation*}
\displaystyle{\large D}_{a}^{\alpha }{\large (f\circ g)(t)=}f^{(\alpha )}%
\big(g(t)\big)\,g^{(\alpha )}(t).
\end{equation*}%
\ \textbf{(v)} All fractional derivatives do not satisfy: $D^{\alpha
}D^{\beta }f=D^{\alpha +\beta }f$, in general \newline
\newline
\textbf{(vi)} All fractional derivatives, specially Caputo definition,
assumes that the function $f$ \ is differentiable.\newline

\bigskip

In \cite{KHKS}, the authors gave a new definition of fractional derivative
which is a natural extension to the usual first derivative as follows:

Given a function $\displaystyle f:[0,\infty )\longrightarrow \mathbb{R}$.
Then for all $t>0,\quad \alpha \in (0,1),$ let
\begin{equation*}
T_{\alpha }(f)(t)=\underset{\varepsilon \rightarrow 0}{\lim }\frac{%
f(t+\varepsilon t^{1-\alpha })-f(t)}{\varepsilon },
\end{equation*}%
$T_{\alpha }$ \ is called \textbf{the conformable fractional derivative of \
}$f$ \textbf{of order }$\alpha .$\newline

Let $f^{(\alpha )}(t)$ stands for $T_{\alpha }(f)(t).$ Hence $\displaystyle %
f^{(\alpha )}(t)=\underset{\varepsilon \rightarrow 0}{\lim }\frac{%
f(t+\varepsilon t^{1-\alpha })-f(t)}{\varepsilon }.$

\textbf{If }$f$\textbf{\ \ is }$\alpha -$\textbf{differentiable in some }$%
(0,b),\;\,b>0$\textbf{, and }$\underset{t\rightarrow 0^{+}}{\lim }f^{(\alpha
)}(t)$ exists$,$ \textbf{then let }
\begin{equation*}
\ f^{(\alpha )}(0)=\underset{t\rightarrow 0^{+}}{\lim }f^{(\alpha )}(t).
\end{equation*}

The conformable derivative satisfies all the classical properties of
derivative.

Further, according to this derivative, the following statements are true,
see \cite{KHKS}.\newline

1. $T_{\alpha }(t^{p})=pt^{p-\alpha }\;$ for all $\;p\in \mathbb{R}$,
\newline

2.\ \ $T_{\alpha }(\sin \frac{1}{\alpha }t^{\alpha })=\cos \frac{1}{\alpha }%
t^{\alpha }$ ,\newline

3.\ \ $T_{\alpha }(\cos \frac{1}{\alpha }t^{\alpha })=-\sin \frac{1}{\alpha }%
t^{\alpha }$ ,\newline

4.\ \ ${\Large T}_{\alpha }{\Large (e}^{\frac{1}{{\Large \alpha }}{\Large t}%
^{{\Large \alpha }}}{\Large )=e}^{\frac{1}{{\Large \alpha }}{\Large t}^{%
{\Large \alpha }}}$.\newline

The\textbf{\ $\alpha-$fractional integral of a function $f$ starting from $%
a\ge 0$}, see \cite{KHKS}, is :
\begin{equation*}
I_{\alpha}^{a}(f)(t)= I_{1}^{a}(t^{\alpha-1}f) = \int_{a}^{t}\frac{f(x)}{%
x^{1-\alpha}}\,dx,
\end{equation*}
where the integral is the usual Riemann improper integral, and $\alpha\in
(0,1)$.\newline
For more about higher conformable fractional integrals and derivatives in left and right senses and other basic concepts we refer to \cite{Th}.

The object of this paper is two folds:\ To introduce the fractional
semigroups of operators associated with the conformable fractional
derivative, then as an application we study the fractional abstract Cauchy
problem according to the conformable fractional derivative which was
introduced in \cite{KHKS}. Indeed, we prove that fractional $\alpha -$
semigroup, is the classical solution of fractional abstract Cauchy problem.
Throughout this paper, \ $\alpha \in (0,1].$

\section{The Basic Definition.}

\begin{definition}
\label{def1} Let $\alpha \in (0,a]$ for any $a>0. $ For a Banach space $X,$
A family $\{T(t)\}_{t\ge 0}\subseteq {\mathcal{L}}(X,X)$ is called a
fractional $\alpha -$semigroup (or $\alpha -$semigroup) of operators if:

\ \ \ \ \ \ $(i)$ $T(0)=I$,

$\ \ \ \ \ \ (ii)$ $T(s+t)^{\frac{1}{\alpha }}=T(s^{\frac{1}{\alpha }})T(t^{%
\frac{1}{\alpha }})$ \ for all $s,t\in \lbrack 0,\infty ).$
\end{definition}

Clearly, if $\alpha =1$, then $1-$semigroups are just the usual semigroups.

\bigskip

\textbf{Example 1.2: }Let $A$ be a bounded linear operator on $X.$ Define $%
T(t)=e^{2\sqrt{t}A}.$ Then $\{T(t)\}_{t\ge 0}$ \ is a $\frac{1}{2}-$%
semigroup. Indeed:

$(i)$ $T(0)=$ $e^{0A}=I.$

$(ii)$ $T(s+t)^{2}=e^{2\sqrt{(s+t)^{2}}%
A}=e^{2(s+t)A}=e^{2sA}e^{2tA}=T(s^{2})T(t^{2}).$

For the definition of conformable fractional exponential matrix and power series expansions see \cite{Th}.
\bigskip

\textbf{Example 2.2. }Let \ $X=C[0,\infty ),$ the space of real valued
continous functions on $[0,\infty ).$ Define $(T(t)f)(s)=f(s+2\sqrt{t}).$
Then one can easily show that $T$ is a $\frac{1}{2}-$ semigroup of operators.

\bigskip

\begin{definition}
\label{def2} An $\alpha -$semigroup $T(t)$ is called a $c_{0}-$semigroup if,
for each fixed $x\in X$, $T(t)x\longrightarrow x \quad$ as $t\longrightarrow
0^{+}$.
\end{definition}

The conformable $\alpha -$derivative of $T(t)$ at $t=0$ \ is called the $%
\alpha -$\textbf{\ infinitesimal generator} of the fractional $\alpha -$%
semigroup $T(t),$ with domain equals
\begin{equation*}
\left\{x\in X: \underset{t\rightarrow 0^{+}}{\lim }T^{(\alpha )}(t)x\;\,%
\mathnormal{exists}\right\}.
\end{equation*}%
We will write $A$ \ for such generator.

\bigskip

\begin{theorem}
\label{thm1} Let $\{T(t)\}_{t\geq 0}\subseteq {\mathcal{L}}(X,X)$ be a $%
c_{0}-\alpha-$semigroup with infinitesimal generator $A$, $0<\alpha \leq 1$.
If $T(t)$ is continuously $\alpha -$differentiable and $x\in D(A)$, then
\begin{equation*}
T^{\alpha }(t)x=AT(t)x=T(t)Ax.
\end{equation*}
\end{theorem}

\textbf{Proof}. Let us begin with
\begin{eqnarray*}
T^{(\alpha )}(t)x&=&\underset{\varepsilon \rightarrow 0}{\lim }\frac{%
T(t+\varepsilon t^{1-\alpha })x-T(t)x}{\varepsilon } \\
&=&\underset{ \varepsilon \rightarrow 0}{\lim }\frac{T(t^{\alpha
}+(t+\varepsilon t^{1-\alpha })^{\alpha }-t^{\alpha })^{\frac{1}{\alpha }%
}x-T(t)x}{ \varepsilon } \\
&=& \underset{\varepsilon \rightarrow 0}{\lim }\frac{T(t^{\alpha
}+((t+\varepsilon t^{1-\alpha })^{\alpha }-t^{\alpha }))^{\frac{1}{\alpha }%
}x-T(t)x}{\varepsilon }\,.
\end{eqnarray*}
Since $T(t)$ is an $\alpha -$semigroup of operators, \ then $T(a+b)^{\frac{1%
}{\alpha }}=T\left( a^{\frac{1}{\alpha }}\right) T\left(b^{\frac{1}{\alpha }%
}\right).$ \ Hence
\begin{eqnarray*}
T^{(\alpha )}(t)x&=&\underset{\varepsilon \rightarrow 0}{ \lim }\frac{%
T(t^{\alpha })^{\frac{1}{\alpha }}T((t+\varepsilon t^{1-\alpha })^{\alpha
}-t^{\alpha }))^{\frac{1}{\alpha }}x-T(t)x}{\varepsilon } \\
&=&\underset{\varepsilon \rightarrow 0}{\lim }\frac{T(t)[T((t+\varepsilon
t^{1-\alpha })^{\alpha }-t^{\alpha }))^{\frac{1}{\alpha }}x-T(0)x]}{%
\varepsilon }.
\end{eqnarray*}
Now, using the \ Mean Value Theorem for conformable fractional derivative,
see \cite{KHKS}, we get
\begin{equation*}
\frac{T(t)[T((t+\varepsilon t^{1-\alpha })^{\alpha }-t^{\alpha }))^{\frac{1
}{\alpha }}x-T(0)x]}{\varepsilon }=T(t)T^{(\alpha )}(c)x\,\frac{%
[(t+\varepsilon t^{1-\alpha })^{\alpha }-t^{\alpha }]}{\alpha \varepsilon }
\end{equation*}
\ for some \ $0<c<(t+\varepsilon t^{1-\alpha })^{\alpha }-t^{\alpha }$.
\newline
\newline
If $\varepsilon \rightarrow 0,$ then $c\rightarrow 0,$ and $\underset{%
\varepsilon \rightarrow 0}{\lim }\ T^{(\alpha )}(c)=T^{(\alpha )}(0)=A.$ \
Consequently,

\begin{equation*}
T^{(\alpha )}(t)x=T(t)Ax\,\underset{\varepsilon \rightarrow 0}{\lim }\frac{%
[(t+\varepsilon t^{1-\alpha })^{\alpha }-t^{\alpha }]}{\alpha \varepsilon }.
\end{equation*}
Using L$^{,}$Hopital's Rule, we get
\begin{equation*}
\underset{\varepsilon \rightarrow 0}{ \lim }\frac{[(t+\varepsilon
t^{1-\alpha })^{\alpha }-t^{\alpha }]}{\alpha \varepsilon }=1.\quad
\mathnormal{\ Hence} \ \ T^{(\alpha )}(t)x=T(t)Ax.
\end{equation*}
Similarly, one can show that $T(t)x\in D(A)\,$ and $\;T^{(\alpha
)}(t)x=AT(t)x $.\newline
This ends the proof.\newline

Let $X=C[0,\infty )$ be the space of continuous real-valued functions such
that $\;\displaystyle \lim_{x\rightarrow \infty}f(x)\;$ is finite, with the
sup norm. Define $T:[0,\infty )\rightarrow {\mathcal{L}}(X,X),$ by
\begin{equation*}
\big(T(t)f\big)(s)=f(s+\frac{1}{\alpha }t^{\alpha }).
\end{equation*}

Claim: \ $T$ \ is an $\alpha -$semigroup.

Indeed:
\begin{eqnarray*}
\big(T(t+k)^{\frac{1}{\alpha }}f\big)(s)&=&f\big(s+\frac{1}{\alpha }[(t+k)^{%
\frac{1}{\alpha } }]^{\alpha }\big) \\
&=&f(s+\frac{1}{\alpha }t+\frac{1}{\alpha }k) \\
&=&\big(T(t^{\frac{1}{\alpha }})T(k^{\frac{1}{ \alpha }})f\big)(s).
\end{eqnarray*}
It is almost immediate that $T(0)=I$ and $T(t)f\in X$ whenever $f\in X$ and
that
\begin{equation*}
\| T(t)f \|_{\infty} \le \| f\|_{\infty}, \quad t\ge 0,
\end{equation*}
so that $T(t)\in {\mathcal{L}}(X,X)$. Since the operator $T(t)$ is a
translation operator corresponding to moving the graph of $\;f\;$ $\frac{1}{%
\alpha}t^{\alpha}$ units to the left and chopping off the part to the left
of the origin, it is known from the literature that $T(t)f$ is
right-continuous at $0$. So $T(t)$ is an $\alpha -$semigroup.\newline

Now,

\begin{theorem}
\label{thm2} The infinitesimal generator of the above semigroup is
\begin{eqnarray*}
&&Af(s)=f^{\prime }(s), \\
&& D(A)=\{ f\in X: f^{\prime}\; \mathnormal{exists \;in }\; X \}.
\end{eqnarray*}
\end{theorem}

\textbf{Proof.}
\begin{eqnarray*}
T^{(\alpha )}(t)f(s)&=&t^{1-\alpha }T^{\prime }(t)f(s) \\
&=&t^{1-\alpha }\underset{\varepsilon \rightarrow 0}{\lim } \frac{%
T(t+\varepsilon )f(s)-T(t)f(s)}{\varepsilon } \\
&=&t^{1-\alpha }\underset{\varepsilon \rightarrow 0}{\lim } \frac{f(s+\frac{1%
}{\alpha }(t+\varepsilon )^{\alpha }-f(s+\frac{1}{\alpha }t^{\alpha })}{%
\varepsilon } \\
&=&t^{1-\alpha }\underset{\varepsilon \rightarrow 0}{\lim } \frac{f(s+\frac{1%
}{\alpha }(t+\varepsilon )^{\alpha })-f(s+\frac{1}{\alpha }t^{\alpha })}{%
\varepsilon }.\,\frac{f(s+t+\varepsilon )-f(s+t)}{f(s+t+\varepsilon )-f(s+t)}
\\
&=&t^{1-\alpha }\underset{\varepsilon \rightarrow 0}{\lim }\frac{f(s+\frac{1%
}{\alpha }(t+\varepsilon )^{\alpha })-f(s+\frac{1}{\alpha }t^{\alpha })}{%
f(s+t+\varepsilon )-f(s+t)}.\,\frac{f(s+t+\varepsilon )-f(s+t)}{\varepsilon }%
\,.
\end{eqnarray*}
Now

$\displaystyle \underset{\varepsilon \rightarrow 0}{\lim } \frac{%
f(s+t+\varepsilon )-f(s+t)}{\varepsilon }=f^{\prime }(s+t)$,

\bigskip $\displaystyle \underset{\varepsilon \rightarrow 0}{\lim } \frac{%
f(s+\frac{1}{\alpha }(t+\varepsilon )^{\alpha }-f(s+\frac{1}{\alpha }%
t^{\alpha })}{f(s+t+\varepsilon )-f(s+t)}=\frac{0}{0}$.\newline

Use L'Hopital's rule (with respect to $\varepsilon )$ to get

\begin{equation*}
\underset{\varepsilon \rightarrow 0}{\lim } \frac{f(s+\frac{1}{\alpha }%
(t+\varepsilon )^{\alpha })-f(s+\frac{1}{ \alpha }t^{\alpha })}{%
f(s+t+\varepsilon )-f(s+t)}=(t)^{\alpha -1}f^{\prime }(s+\frac{1}{\alpha }%
(t)^{\alpha })\big/ f^{\prime }(s+t).
\end{equation*}
Thus the product gives
\begin{equation*}
T^{(\alpha )}(t)f(s)=f^{\prime }(s+\frac{1}{\alpha }t^{\alpha}).
\end{equation*}
Now take the limit as $t\rightarrow 0$ \ to get
\begin{equation*}
T^{(\alpha )}(0)f(s)=f^{\prime }(s).
\end{equation*}
Hence \ \ $Af=f$ $^{\prime }$. This completes the proof.\newline

Let us show how our theory can be applied to obtain information about
solutions of certain problems. In particular, we want to use the fractional
semigroups approach for solving the so-called $\alpha-$ abstract Cauchy
problem.

\begin{definition}
\label{def3} Let $X$ be a Banach space, $A:D(A)\subseteq X\longrightarrow X$
a linear operator and $u_{0}\in X$. A function $\;u:[0,\infty)%
\longrightarrow X\; $ is a solution of the $\alpha-$ abstract Cauchy problem
\begin{eqnarray}
&&u^{(\alpha)}(t)=A u(t), \quad t>0,  \label{abc1} \\
&&u(0)=u_{0}  \label{abc2}
\end{eqnarray}
if:

\ \ \ \ \ \ $(i)$ $u$ is continuous on $[0,\infty)$,

$\ \ \ \ \ \ (ii)$ $u$ is continuously $\alpha-$differentiable on $%
(0,\infty) $,

$\ \ \ \ \ \ (iii)$ $u(t)\in D(A)\quad$ for $\; t>0$,

$\ \ \ \ \ \ (iv)$ $u$ satisfies (\ref{abc1})-(\ref{abc2}).
\end{definition}

We have

\begin{theorem}
\label{thm3} Let $X$ be a Banach space and $A$ the infinitesimal generator
of a $c_{0}-\alpha-$semigroup $\{T(t)\}_{t\ge 0}\subseteq {\mathcal{L}}(X,X)$%
. If $u_{0}\in D(A)$, then problem (\ref{abc1})-(\ref{abc2}) has one and
only one solution $u$, namely,
\begin{equation*}
u(t)=T(t)u_{0}.
\end{equation*}
\end{theorem}

\textbf{Proof.} Clearly $u(t)=T(t)x$ is a solution of problem (\ref{abc1})-(%
\ref{abc2}). For uniqueness, let $u$ be a solution of (\ref{abc1})-(\ref%
{abc2}). Then
\begin{eqnarray*}
[T(t-s)u(s)]^{(\alpha)}&=& T(t-s)u^{(\alpha)}(s)-AT(t-s)u(s) \\
&=& T(t-s)u^{(\alpha)}(s)-T(t-s)Au(s) \\
&=& T(t-s)[u^{(\alpha)}(s)-Au(s)] \\
&=& 0.
\end{eqnarray*}
Applying $I_{\alpha}^{0}$, in $s$, gives
\begin{equation*}
T(t-t)u(t)-T(t)u_{0}=0\Rightarrow u(t)=T(t)u_{0}.
\end{equation*}
This completes the proof.

\begin{remark}
\label{rem1} Let $X=C[0,\infty )$ be the space of continuous real-valued
functions such that $\;\displaystyle \lim_{x\rightarrow \infty}f(x)\;$ is
finite, with the sup norm. Define the operator $A$ by
\begin{eqnarray*}
&&Af(s)=f^{\prime }(s), \\
&& D(A)=\{ f\in X: f^{\prime}\; \mathnormal{exists\; in }\; X \}.
\end{eqnarray*}
Then $A$ is a generator of the above $\alpha-$semigroup of translation. If $%
u_{0}\in D(A)$, then problem (\ref{abc1})-(\ref{abc2}) has the unique
solution $u(t)=T(t)u_{0}$, where $T(t)$ is an $\alpha-$semigroup generated
by $A$.

Now it is readily seen that if $g$ is continuously differentiable on $%
[0,\infty)$, then
\begin{equation*}
u(x,t)=g\left( x+\frac{1}{\alpha}t^{\alpha} \right)
\end{equation*}
is the unique solution of the problem
\begin{eqnarray*}
&&\frac{\partial^{\alpha} u}{ \partial t^{\alpha}}=\frac{\partial u}{
\partial x},\quad x>0,\;\; t>0, \\
&& u(x,0)=g(x),\quad x>0.
\end{eqnarray*}
\end{remark}


\begin{thebibliography}{9}
\bibitem{Th} Abdeljawad, T., On Conformable Fractional Calculus,Journal of Computational ad Applied Mathematics,Vol. 279, 1 May 2015, 57-66, arXiv: 1402.6892v1 [math D, S] 27 Feb 2014.
\bibitem{H2010} Al Horani, M., Projection Method for Solving Degenerate
First-Order Identification Problem. J. Math. Anal. Appl. Vol. 364, pp.
204-208, 2010.

\bibitem{Bala} Balakrishnan, A. V., \ Fractional Powers Of Closed Operators
And The Semigroups Generated By Them, Pacific Journal of Mathematics 10, pp.
419-439, 1960.

\bibitem{Boch} Bochner, S., Diffusion Equation And Stochastic Processes,
Proc. Nat. Acad. Sciences, U.S.A., 35, pp. 368-370, 1949.

\bibitem{EN} Engel, K.-J. and Nagel, R:\textit{\ One Parameter Semigroups
For Linear Evolution Equations}, Graduate Texts in Math., Springer-Verlag,
Berlin-Heidelgerg-New York, 2000.

\bibitem{KHKS} Khalil, R., Al Horani, M., Yousef. A. and Sababheh, M., A new
Definition Of Fractional Derivative, J. Comput. Appl. Math. 264. pp. 65--70,
2014.

\bibitem{Paz} Pazy, A., \textit{Semigroups of Linear Operators and
Applications to Partial Differntial Equations}, Springer-Verlag, 1983.

\bibitem{Pop} Popescu, E., On The Fractional Cauchy Problem Associated With
a Feller Semigroup, Romanian Academy, Mathematical Reports, Vol. 12(62), No.
2, pp. 181--188, 2010.

\bibitem{POV} Prilepko, I., Orlovsky, G. and Vasin, A., \textit{Methods For
Solving Inverse Problems In Mathematical Physics}. Marcel Dekker. Inc. New
York, 2002.
\end{thebibliography}
\end{document}